\documentclass[
final
]{dmtcs-episciences}

\usepackage[utf8]{inputenc}
\usepackage{subfigure}

\usepackage[round]{natbib}

\usepackage{amssymb}
\usepackage{amssymb}
\usepackage{enumerate}
\usepackage{latexsym}
\usepackage[only,ninrm,elvrm,twlrm,sixrm,egtrm,tenrm]{rawfonts}
\usepackage{indentfirst}
\usepackage{amsmath,amsfonts}
\usepackage{algorithm}
\usepackage{algorithmic}
\usepackage{multirow}
\usepackage{graphicx,psfrag}
\usepackage{graphics}
\usepackage{makeidx}
\usepackage{ifpdf}
\usepackage{amsmath}
\usepackage{xcolor}
\usepackage{array}
\usepackage[all,cmtip,color]{xypic}
\usepackage{lipsum}

\newtheorem{thm}{Theorem}
\newtheorem{lemma}{Lemma}
\newtheorem{defn}{Definition}

\def\qed{\hfill \nopagebreak\rule{5pt}{8pt}}

\author{Zongwen Bai\affiliationmark{1,2}\thanks{Supported by the Natural Science Foundation of Shaanxi province of China (No. 2014JM8357, 2014JQ2-6031).}
  \and Jianhua Tu\affiliationmark{3}\thanks{Corresponding author.}
  \and Yongtang Shi\affiliationmark{4}\thanks{Supported by National Natural Science Foundation of China,
Natural Science Foundation of Tianjin (No. 17JCQNJC00300), the China-Slovenia bilateral project ``Some topics in modern graph theory" (No.~12-6),
Open Project Foundation of Intelligent Information Processing Key Laboratory of Shanxi Province (No. CICIP2018005),
and the Fundamental Research Funds for the Central Universities.}}
\title{An improved algorithm for the vertex cover $P_3$ problem on graphs of bounded treewidth}
\affiliation{
  School of Computer Science, Northwestern Polytechnical University, Xi'an, Shaanxi, 710072, China\\
  School of Physics and Electronic Information, Yan'an University, Yan'an, Shaanxi 716000, China\\
  Department of Mathematics, Beijing University of Chemical Technology, Beijing 100029, China\\
  Center for Combinatorics and LPMC, Nankai University, Tianjin 300071, China}
\keywords{Combinatorial optimization, Vertex cover $P_3$ problem, Connected vertex cover $P_3$ problem, Treewidth, Dynamic programming}
\received{2016-4-2}
\revised{2017-6-30, 2018-12-27, 2019-10-11}
\accepted{2019-10-11}
\begin{document}
\publicationdetails{21}{2019}{4}{17}{1425}
\maketitle
\begin{abstract}
Given a graph $G=(V,E)$ and a positive integer $t\geq2$, the task in the vertex cover $P_t$ ($VCP_t$) problem is to find a minimum subset of vertices $F\subseteq V$ such that every path of order $t$ in $G$ contains at least one vertex from $F$. The $VCP_t$ problem is NP-complete for any integer $t\geq2$.  Recently, the authors presented a dynamic programming algorithm with runtime $4^{p}\cdot n^{O(1)}$ that can solve the $VCP_3$ problem in any $n$-vertex graph given together with its tree decomposition of width at most $p$. In this paper, we propose an improvement of it and improved the time-complexity to $3^{p}\cdot n^{O(1)}$.

The connected vertex cover $P_3$ ($CVCP_3$) problem is the connected variation of the $VCP_3$ problem where $G[F]$ is required to be connected. Using the Cut\&Count technique, we give a randomized algorithm with runtime $4^{p}\cdot n^{O(1)}$ that can solve the $CVCP_3$ problem in any $n$-vertex graph given together with its tree decomposition of width at most $p$.
\end{abstract}

\section{Introduction}

In this paper, we consider only finite, simple and undirected graphs $G$. For a graph $G$, we denote its vertex set and edge set by $V(G)$ and $E(G)$, respectively. Unless stated otherwise, let $n:=|V(G)|$ and $m:=|E(G)|$. As usual, $P_t$ denotes the path on $t$ vertices. Given a graph $G$ and a positive integer $t\geq2$, a subset of vertices $F$ is called a vertex cover $P_t$ ($VCP_t$) set if every path of order $t$ contains at least one vertex from $F$. The task in the $VCP_t$ problem is to find a minimum $VCP_t$ set in $G$.

The $VCP_t$ problem is NP-complete for any integer $t\geq2$ \cite{lewis}. Clearly, the $VCP_2$ problem is the well-known vertex cover problem. Thus, the $VCP_t$ problem is a natural extension of the vertex cover problem.

In this paper, we restrict our attention to the $VCP_3$ problem. In the literature, the $VCP_3$ problem is also known as the
Bounded-Degree-one Deletion (1-BDD) problem. Given a graph $G=(V,E)$, the task in the 1-BDD problem is to find a minimum subset of vertices $F$ such that the graph $G[V\setminus F]$ has maximum degree at most one. For the $VCP_3$ problem, Kardo\v{s} et al. \cite{kar} gave an exact algorithm with runtime $1.5171^n\cdot n^{O(1)}$ and a randomized approximation algorithm with an expected approximation ratio 23/11. Chang et al. \cite{chang} presented an improved exact algorithm with runtime $1.4658^n \cdot n^{O(1)}$. For weighted version of the $VCP_3$ problem, Tu and Zhou \cite{tu2,tu3} presented two 2-approximation algorithms using the primal-dual method and the local-ratio method.

On the other hand, many efforts have also been made in study of parameterized complexity of the $VCP_3$ problem with the size $k$ of the solution as the parameter. Using the iterative compression technique, a parameterized algorithm running in time $2^k\cdot n^{O(1)}$ was given independently by Fellows et al. \cite{fel} and Moser et al. \cite{moser}. Wu \cite{wu} gave an improved algorithm with runtime $1.882^k\cdot n^{O(1)}$ using the Measure \& Conquer approach. Katreni\v{c} \cite{kat} gave further
improvement so that running time reaches $1.8172^k\cdot n^{O(1)}$.

In this paper, we consider the $VCP_3$ problem on graphs of bounded treewidth. The treewidth ($tw$) is a graph parameter that plays a fundamental
role in various graph algorithms. Very roughly, treewidth captures how similar a graph is to a tree. It is well-known that many NP-hard graph problem can be solved efficiently if the input graph $G$ has bounded treewidth. For example, an algorithm to solve the vertex cover problem running in time $4^{tw}\cdot n^{O(1)}$ was given in \cite{klein}, while the book \cite{nie} presented an algorithm running in time $2^{tw}\cdot n^{O(1)}$. For one more example, Alber et al. \cite{alber2} gave a $4^{tw}\cdot n^{O(1)}$ time algorithm for the dominating set problem, improving
over the natural $9^{tw}\cdot n^{O(1)}$ algorithm of Telle and Proskurowski \cite{telle}. Using fast subset convolution \cite{bjo,cygan}, the running time of the algorithm of Alber et al. can be improved and reaches $3^{tw}\cdot n^{O(1)}$.

Recently, Tu et al. \cite{tu4} presented a dynamic programming algorithm running in time $4^{tw}\cdot n^{O(1)}$ for the $VCP_3$ problem on graphs of bounded treewidth.

\begin{thm}\cite{tu4}
Let $G$ be an $n$-vertex graph given together with its tree decomposition of width at most $p$. Then the $VCP_3$ problem in $G$ can be solved in time $4^p\cdot n^{O(1)}$.
\end{thm}

Some families of graphs with bounded treewidth include the cactus graphs, pseudoforests, series-parallel graphs, outerplanar graphs, Halin graphs,  Apollonian networks, and so on. Thus, the $VCP_3$ problem on these families of graphs can be solved in polynomial time. Moreover, in order to solve the $VCP_3$ problem on these families of graphs more efficiently, it is necessary to obtain an improved algorithm on graphs of bounded treewidth.

We improve the result in \cite{tu4} in the following ways.

\medskip
1. Using a refined variant of nice tree decompositions, we present a new dynamic programming algorithm running in time $3^{tw}\cdot n^{O(1)}$ for the $VCP_3$ problem.

\medskip
2. Consider the connected vertex cover $P_3$ ($CVCP_3$) problem, a variation of the $VCP_3$ problem where it is required that the graph
induced by the $VCP_3$ set is connected. Using the Cut\&Count technique, we give a randomized algorithm with runtime $4^{tw}\cdot n^{O(1)}$ for the $CVCP_3$ problem.
\medskip

It should be pointed out that the dynamic programming and the Cut\&Count technique are two widely used methods. The dynamic programming can give exact algorithms running in time $2^{O(tw)}\cdot n^{O(1)}$ for many NP-hard problems on graphs with a given tree decomposition of width $tw$ \cite{cygan}, while the Cut\&Count method can give randomized algorithms for a wide range of connectivity problems running in time $c^{O(tw)}\cdot n^{O(1)}$ for a small constant $c$ \cite{cygan2}. The main work of this paper is to give related results for the $VCP_3$ problem and its connected variation using the two methods.

In \cite{pil}, Pilipczuk introduced the existential counting modal logic (ECML), then extended ECML by adding connectivity requirements and proved that problems expressible in the extension are tractable in single exponential time when parameterized by treewidth; however, using randomization. Thus, the work by Pilipczuk
already shows that $CVCP_3$ can be solved in single-exponential time. Our result provides a better exponential dependence than the very general result in \cite{pil}.

The remaining part of this paper is organized as follows.
In Section 2, we give some notation and introduce formally the concepts of tree decomposition and treewidth.
In Section 3, we present a new dynamic programming algorithm for the $VCP_3$ problem on graphs of bounded treewidth. In Section 4, a randomized algorithm for the $CVCP_3$ problem is given.

\section{Preliminaries}

In this section we give some notation which we make use of in the paper. Let $G=(V,E)$ be a graph. For $v\in V$, denote by $N_G(v)$ the set of neighbors of $v$ in $G$. Let $d_G(v):=|N_G(v)|$. For a subset $V'\subseteq V$, we denote by $G[V']$ the subgraph induced by $V'$ and write $G-V'$ as an abbreviation for the induced subgraph $G[V\setminus V']$. For all terminology and notation not defined here, we refer the reader to \cite{bon}.

\begin{defn}
A {\it tree decomposition} of a graph $G=(V,E)$ is a pair $\mathcal{T}=(T,\{X_i:i\in I\})$, where each $X_i\subseteq V$ is called a bag, and $T$ is a rooted tree with the elements of $I$ as nodes. The following three conditions must hold:
\begin{enumerate}
\item[1.] $\bigcup_{i\in I}X_i=V$,
\item[2.] for all edges $uv\in E$ there exists an $i\in I$ with $u\in X_i$ and $v\in X_i$, and
\item[3.] for all $i,j,k\in I$, if $j$ lies on the path from $i$ to $k$ in $T$ then $X_i\cap X_k\subseteq X_j$.
\end{enumerate}
\end{defn}
The {\it width} of a tree decomposition $\mathcal{T}=(T, \{X_i:i\in I\})$ is $\max_{i\in I}|X_i|-1$. The {\it treewidth} $tw(G)$ of a graph $G$ is defined as the minimum width over all tree decompositions of $G$. A tree has treewidth one, many well-studied graph families also have bounded treewidth.

Given a graph $G=(V,E)$, if $p$ is a fixed constant, the problem to determine whether the treewidth of $G$ is at most $p$ can be decided in linear time and a corresponding tree decomposition can be constructed in linear time with a high constant factor \cite{bod1}. There are also several good heuristic algorithms which can construct tree decompositions of small width, if existing, and often work well in practice (see \cite{bod3}).

Given an $n$-vertex graph $G=(V,E)$, a tree decomposition $\mathcal{T}=(T, \{X_i:i\in I\})$ can be converted (in polynomial time) in a {\it nice} tree decomposition of the same width $p$ and with $O(pn)$ nodes \cite{klo}: here the tree $T$ is rooted and binary, and the following conditions are satisfied.
\begin{itemize}
\item  All the leaves as well as the root of $T$ contain empty bags.
\item  All non-leaf nodes of $T$ are of three types:
\begin{itemize}
\item {\bf Introduce node:} a node $t$ with exactly one child $t'$ such that $|X_t|=|X_{t'}|+1$ and $X_t=X_{t'}\cup\{v\}$ for some vertex $v\notin X_{t'}$.
\item {\bf Forget node:} a node $t$ with exactly one child $t'$ such that $|X_t|=|X_{t'}|-1$ and $X_t=X_{t'}\setminus\{v\}$ for some vertex $v\in X_{t'}$.
\item {\bf Join node:} a node $t$ with two children $t_1$ and $t_2$ such that $X_t=X_{t_1}=X_{t_2}$.
\end{itemize}
\end{itemize}

\section{An improved algorithm for the $VCP_3$ problem on graphs of bounded treewidth}

In this section, we will use a refined variant of nice tree decompositions to obtain an improved dynamic programming algorithm for solving the $VCP_3$ problem in graphs of bounded treewidth. Given a graph $G=(V,E)$ and its nice tree decomposition $\mathcal{T}=(T,\{X_i:i\in I\})$, we will add a new type of a node called an {\it introduce edge node} in the refined variant of nice decomposition, which is defined as follows.
\begin{itemize}
\item {\bf Introduce edge node:} an internal node $t$ of $T$, labeled with an edge $uv\in E(G)$ with exactly one child $t'$ for which $u,v\in X_t=X_{t'}$. We say that edge $uv$ is introduced at $t$.
\end{itemize}

We additionally require that every edge in $E$ is introduced exactly once. The introduce edge node enables us to add edges one by one. A standard nice tree decomposition of width $p$ can be transformed to the refined variant with the same width $p$ and with $O(pn)$ nodes in polynomial time \cite{cygan}.

Let $\mathcal{T}=(T,\{X_i:i\in I\})$ be a refined variant of nice tree decomposition. Recall that then $T$ is rooted at some node $r$. For each node $t$ of $T$, we denote by $V_t$ the union of all the bags present in the subtree of $T$ rooted at $t$.  Associate a subgraph $G_t$ of $G$ with each node $t$ of $T$ defined as follows.
$$G_t=(V_t,E_t=\{e: e \text{ is introduced in the subtree rooted at t}\})$$

Now, we present an improved dynamic programming algorithm on the refined variant of nice tree decomposition for solving the $VCP_3$ problem and give the following theorem.

\begin{thm}
Let $G$ be an $n$-vertex graph given together with its refined variant of nice tree decomposition of width at most $p$. Then the $VCP_3$ problem in $G$ can be solved in time $O(3^{p}\cdot pn^{6})$.
\end{thm}

\begin{proof}
Let $\mathcal{T}=(T,\{X_i\}_{i\in I})$ be a refined variant of nice tree decomposition of $G$. Assume that $r$ is the root node of $T$.

For each node $t$ of $T$, a coloring of bag $X_t$ is a mapping $f\colon X_t\rightarrow\{1,0_0,0_1\}$ assigning three different colors to vertices of the bag. Clearly, there exist $3^{|X_t|}$ colorings of $X_t$. For a coloring $f$ of $X_t$, denote by $c[t,f]$ the minimum size of a $VCP_3$ set $F\subseteq V_t$ in $G_t$ so that
\begin{itemize}
\item $F\cap X_t=f^{-1}(1)$. The meaning is that all vertices of $X_t$ colored 1 have to be contained in $F$,
\item each vertex colored $0_0$ is an isolated vertex in $G_t-F$,
\item each vertex colored $0_1$ has degree 1 in $G_t-F$.
\end{itemize}
We put $c[t,f]=+\infty$ if no such $VCP_3$ set $F$ for $t$ and $f$ exists. Note that because $G_r=G$ and $X_r=\emptyset$, $c[r,\emptyset]$ is exactly the minimum size of a $VCP_3$ set in $G$.

Now, we give the recursive formulas for the values of $c[t,f]$. For each node $t$ and a coloring $f$ of $X_t$, we compute the value of $c[t,f]$ based on the values computed for the children of $t$. By applying the formulas in a bottom-up manner on $T$, we can obtain the value of $c[r,\emptyset]$, which is exactly the minimum size of a $VCP_3$ set in $G$.

\medskip
\noindent {\bf Leaf node.}

For a leaf node $t$, $X_t=\emptyset$. Hence, there is only one empty coloring for $X_t$ and $c[t,\emptyset]=0$.

\medskip
\noindent{\bf Introduce node.}

Let $t$ be an introduce node with a child $t'$ so that $X_t=X_{t'}\cup\{v\}$ for some $v\notin X_{t'}$. Note that $V_t=V_{t'}\cup\{v\}$ and $E_t=E_{t'}$, i.e., $v$ is an isolated vertex in $G_t$. Hence, for a coloring $f$ of $X_t$, we just need to be sure that if $F$ is a $VCP_3$ set for $f$ and $t$ and $v$ is not contained in $F$, then $v$ must be an isolated vertex in $G_t-F$.

For a coloring $f$ of $X_t$, we denote by $f|_{X_{t'}}$ the restriction of $f$ to $X_{t'}$. It is easy to see that the following formula holds.

$$c[t,f]=\left\{
         \begin{array}{ll}
           c[t',f|_{X_{t'}}]+1, & \hbox{when $f(v)=1$,}\\
           c[t',f|_{X_{t'}}] & \hbox{when $f(v)=0_0$,}\\
           +\infty, & \hbox{when $f(v)=0_1$}.
         \end{array}\right.$$

\medskip

\noindent{\bf Introduce edge node.}

Let $t$ be an introduce edge node labeled with an edge $uv\in E(G)$ and $t'$ be the child of $t$. Note that $V_t=V_{t'}$, $E_{t}=E_{t'}\cup\{uv\}$ and $uv\notin E_{t'}$. In other words, the only difference between the graph $G_t$ and the graph $G_{t'}$ is that $uv\in E(G_t)$ and $uv\notin E(G_{t'})$. Hence, for a coloring $f$ of $X_t$, we just need to be sure that if $F$ is a $VCP_3$ set for $f$ and $t$ and neither $u$ nor $v$ is contained in $F$, then $u$ and $v$ are vertices of degree 1 in $G_t-F$. Thus, for a coloring $f$ we can obtain the following formula.
$$c[t,f]=\left\{
         \begin{array}{ll}
           c[t',f], & \hbox{when $f(u)=1$ or $f(v)=1$,}\\
           c[t',f'] & \hbox{when $(f(u),f(v))=(0_1,0_1)$,}\\
           +\infty, & \hbox{otherwise},
         \end{array}\right.$$
where $f'$ is a coloring of $X_{t'}$ and $f'(x)=\left\{\begin{array}{ll}f(x), &\hbox{when $x\neq u$ and $x\neq v$,}\\0_0, &\hbox{when $x=u$ or $x=v$.}\end{array}\right.$

\medskip
\noindent{\bf Forget node.}

Let $t$ be a forget node with a child $t'$ so that $X_t=X_{t'}\setminus\{v\}$ for some $v\in X_{t'}$. Note that $V_t=V_{t'}$ and $G_t=G_{t'}$. For a coloring $f$ of $X_t$ and a color $\alpha\in\{1,0_0,0_1\}$, we define a coloring $f_{v\rightarrow \alpha}$ of $X_{t'}$ as follows:

$$f_{v\rightarrow\alpha}(x)=\left\{
         \begin{array}{ll}
           f(x), & \hbox{when $x\neq v$,} \\
           \alpha, & \hbox{when $x=v$.}
         \end{array}
       \right.
$$

We claim that the following formula holds.

$$c[t,f]=\min\{c[t',f_{v\rightarrow 1}],c[t',f_{v\rightarrow 0_0}],c[t',f_{v\rightarrow 0_1}]\}$$

Suppose that $F$ is a $VCP_3$ set for which the minimum is attained in the definition of $c[t,f]$. Thus, for any vertex $u\in V_t$, either $u\in F$ or $d_{G_t-F}(u)\leq 1$. If $v\in F$, then $F$ is one of the sets considered in the definition of $c[t',f_{v\rightarrow 1}]$, and hence $c[t',f_{v\rightarrow 1}]\leq c[t,f]$. For $k\in\{0,1\}$, if $v\notin F$ and $v$ is a vertex of degree $k$ in $G_t-F$, then $F$ is one of the sets considered in the definition of $c[t',f_{v\rightarrow 0_k}]$, and hence $c[t',f_{v\rightarrow 0_k}]\leq c[t,f]$. Thus
$\min\{c[t',f_{v\rightarrow 1}],c[t',f_{v\rightarrow 0_0}],c[t',f_{v\rightarrow 0_1}]\}\leq c[t,f]$.

On the other hand, each set that is considered in the definition of $c[t',f_{v\rightarrow 1}]$ is also considered in the definition of $c[t,f]$, and the same holds also for $c[t',f_{v\rightarrow 0_0}]$ and $c[t',f_{v\rightarrow 0_1}]$. This means $c[t,f]\leq c[t',f_{v\rightarrow 1}]$, $c[t,f]\leq c[t',f_{v\rightarrow 0_0}]$ and $c[t,f]\leq c[t',f_{v\rightarrow 0_1}]\}$.

In conclusion, $c[t,f]=\min\{c[t',f_{v\rightarrow 1}],c[t',f_{v\rightarrow 0_0}],c[t',f_{v\rightarrow 0_1}]\}$.

\medskip
\noindent{\bf Join node.}

Suppose that $t$ is a join node with children $t_1$, $t_2$ so that $X_t=X_{t_1}=X_{t_2}$. Note that $V_t=V_{t_1}\cup V_{t_2}$, $V_{t_1}\cap V_{t_2}=X_t$ and $E_t=E_{t_1}\cup E_{t_2}$. It is important to point out that because every edge in $E(G)$ is introduced exactly once in the refined variant of nice tree decomposition, $E_{t_1}\cap E_{t_2}=\emptyset$.

We say that a pair of colorings $f_1$ of $X_{t_1}$, $f_2$ of $X_{t_2}$ is consistent with a coloring $f$ of $X_t$ if for each vertex $v\in X_t$, the following conditions hold.
\begin{description}
  \item[(1)] $f(v)=1$ if and only if $f_1(v)=f_2(v)=1$,
  \item[(2)] $f(v)=0_0$ if and only if $f_1(v)=f_2(v)=0_0$,
  \item[(3)] $f(v)=0_{1}$ if and only if $(f_1(v),f_2(v))\in\{(0_{1},0_0),(0_0,0_{1})\}$.
\end{description}

We will show that the following recursive formula holds:
$$c[t,f]=\min_{f_1,f_2}\{c[t_1,f_1]+c[t_2,f_2]\}-|f^{-1}(1)|,$$
where the minimum is taken over all pairs of colorings $f_1$ (of $X_{t_1}$), $f_2$ (of $X_{t_2}$) consistent with $f$.

On one hand, if $F$ is a $VCP_3$ set for $t$ and $f$, then $F\cap V_{t_1}$ and $F\cap V_{t_2}$ are $VCP_3$ sets for $t_1$ and $f_1$, and $t_2$ and $f_2$, for some pairs of colorings $f_1,f_2$ that are consistent with $f$. The reason is as follows: For a vertex $v\in X_t$, if $v\in F\cap X_t$, then $v\in F\cap X_{t_1}$ and $v\in F\cap X_{t_2}$; if $v\notin F$ and $d_{G_t-F}(v)=0$, then $d_{G_{t_1}-F\cap V_{t_1}}(v)=0$ and $d_{G_{t_2}-F\cap V_{t_2}}(v)=0$; if $v\notin F$ and $d_{G_t-F}(v)=1$, then $v$ is a 1-degree vertex in $G_{t_1}-F\cap V_{t_1}$ and an isolated vertex in $G_{t_2}-F\cap V_{t_2}$, or $v$ is an isolated vertex in $G_{t_1}-F\cap V_{t_1}$ and a 1-degree vertex in $G_{t_2}-F\cap V_{t_2}$, because there exist no edges between $V_{t_1}\setminus X_t$ and $V_{t_2}\setminus X_t$, $d_{G_t-F}(v)\leq 1$ for any vertex $v\in V_t\setminus (X_t\cup F)$. Hence,
\begin{center}
$\min_{f_1,f_2}\{c[t_1,f_1]+c[t_2,f_2]\}\leq |F\cap V_{t_1}|+|F\cap V_{t_2}|=|F|+|F\cap X_t|=c[t,f]+|f^{-1}(1)|,$
\end{center}
where the minimum is taken over all pairs of colorings $f_1$, $f_2$ consistent with $f$.

On the other hand, given a coloring $f$ of $X_t$, consider any pair of colorings $f_1,f_2$ that is consistent with $f$. Suppose that $F_1$ and $F_2$ are $VCP_3$ sets for $t_1$ and $f_1$, $t_2$ and $f_2$. Let $F:=F_1\cup F_2$. Since $f_1^{-1}(1)=f_2^{-1}(1)=f^{-1}(1)$, $F_1\cap X_{t_1}=F_2\cap X_{t_2}=F\cap X_t$ and $|F|=|F_1|+|F_2|-f^{-1}(1)$. Recall that the pair of colorings $f_1$, $f_2$ is consistent with $f$. For every vertex $v\in X_t\setminus F$, $d_{G_t-F}(v)\leq 1$, because there exist no edges between $V_{t_1}\setminus X_t$ and $V_{t_2}\setminus X_t$, $d_{G_t-F}(v)\leq 1$ for any vertex $v\in V_t\setminus (X_t\cup F)$. Thus, $F$ is a $VCP_3$ set for $t$ and $f$. Hence, $c[t,f]\leq |F|=|F_1|+|F_2|-f^{-1}(1)$. Since the pair of colorings can be chosen arbitrarily,
$c[t,f]\leq \min_{f_1,f_2}\{c[t_1,f_1]+c[t_2,f_2]\}-|f^{-1}(1)|,$ where the minimum is taken over all pairs of colorings $f_1$, $f_2$ consistent with $f$. We thus have proved the above recursive formula.

Now, we have finished the description of the recursive formulas for the values of $c[t,f]$, i.e., we have presented an improved dynamic programming algorithm for the $VCP_3$ problem on graphs of bounded treewidth. The value of $c[t,\emptyset]$ is exactly the minimum size of a $VCP_3$ set in $G$. Moreover, when backtracking how the value of $c[t,\emptyset]$ is obtained, we can construct a minimum $VCP_3$ set $F$ in $G$.

Next, we analyze the running time of the algorithm.
Clearly, for each leaf node, introduce node, introduce edge node, forget node, the time needed to computer all values of $c[t,f]$ is $3^p$. We will apply fast algorithms for subset convolution to compute the time needed to process each join node.

Given a set $S$ and two functions $g,h\colon 2^S\rightarrow \mathbb{Z}$, the {\it subset convolution} of $g$ and $h$ is a function
$(g\ast h)\colon 2^S\rightarrow \mathbb{Z}$ so that for every $Y\subseteq S$,
$$(g\ast h)(Y)=\min_{\mbox {\tiny $\begin{array}{c}A\cup B=Y\\A\cap B=\emptyset\end{array}$}}(g(A)+h(B)).$$

\begin{lemma}\label{lemma1}\cite{bjo,van}
Let $S$ be a set with $n$ elements and $M$ be a positive integer.
For two functions $g,h\colon 2^S\rightarrow\{-M,\cdots,M\}\cup\{+\infty\}$, if all the
values of $g$ and $h$ are given, then all the $2^n$ values of the subset convolution of
$g$ and $h$ can be computed in $2^n\cdot n^3\cdot O(M\log(Mn)\log\log(Mn))$ time.
\end{lemma}

Let $t$ be a join node with the children $t_1$ and $t_2$. For a coloring $f$ of $X_t$, if colorings $f_1$ of $X_{t_1}$ and $f_2$ of $X_{t_2}$ are consistent with $f$, then the following conditions hold.

\begin{itemize}
  \item $f^{-1}(1)=f_1^{-1}(1)=f_2^{-1}(1)$,
  \item $f^{-1}(0_1)=f_1^{-1}(0_1)\cup f_2^{-1}(0_1)$,
  \item $f_1^{-1}(0_1)\cap f_2^{-1}(0_1)=\emptyset$.
\end{itemize}

Fix a set $R\subseteq X_t$, let $\mathcal{F}_R$ denote the set of all functions $f\colon X_t\rightarrow\{1,0_0,0_1\}$ so that $f^{-1}(1)=R$. Observe that every coloring $f\in \mathcal{F}_R$ can be represented by a set $Y\subseteq X_t\setminus R$ as follows.

$$g_Y(x)=\left\{
           \begin{array}{ll}
             1, & \hbox{if $x\in R$;} \\
             0_1, & \hbox{if $x\in Y$;} \\
             0_0, & \hbox{if $x\in X_t\setminus(R\cup Y)$.}
           \end{array}
         \right.
$$
Then for every $f(=g_Y)\in\mathcal{F}_{R}$, the following formula holds.

$$c[t,f]=c[t,g_Y]=\min_{\mbox {\tiny $\begin{array}{c}A\cup B=Y\\A\cap B=\emptyset\end{array}$}}(c[t_1,g_A]+c[t_2,g_B])-|R|.$$

Recall that $n$ is the number of vertices in $G$. We can define two functions $h_1,h_2\colon 2^{X_{t}\setminus R}\rightarrow\{0,\cdots,n\}\cup \{+\infty\}$ so that for every set $T\subseteq X_t\setminus R$ we have $h_1(T)=c[t_1,g_T]$ and $h_2(T)=c[t_2,g_T]$. Thus,

$$c[t,f]=c[t,g_Y]=\min_{\mbox {\tiny $\begin{array}{c}A\cup B=Y\\A\cap B=\emptyset\end{array}$}}(h_1(A)+h_2(B))-|R|=(h_1\ast h_2)(Y)-|R|,$$
where $h_1\ast h_2$ is the subset convolution of $h_1$ and $h_2$. By Lemma \ref{lemma1}, we can compute all the $2^{|X_{t}\setminus R|}$ values of $c[t,f]$ for all $f\in \mathcal{F}_{R}$ in time $2^{|X_t\setminus R|}\cdot O(n^{5})$. Since $\sum_{R\subseteq X_t}2^{|X_t\setminus R|}=3^{|X_t|}\leq 3^{p+1}$, the time needed to process each joint node is $O(3^p\cdot n^{5})$.

Recall that the number of nodes of $T$ in the tree decomposition is $O(pn)$, the running time of the algorithm presented here
is $O(3^{p}\cdot pn^{6})$.
\end{proof}

\section{A randomized algorithm for the connected vertex cover $P_3$ problem}

In this section we consider the connected vertex cover $P_3$ ($CVCP_3$) problem, a variation of the $VCP_3$
problem where it is required that the graph induced by the $VCP_3$ set is connected. We present a randomized algorithm with runtime $4^{tw}\cdot n^{O(1)}$ for the $CVCP_3$ problem on $n$-vertex graphs with treewidth $tw$.
The main method used here is the Cut\&Count technique. The technique was introduced by Cygan et al. \cite{cygan2} and can give randomized algorithms for many  problems with certain connectivity requirement running in time $c^{tw}n^{O(1)}$ for a small constant $c$. We show that for the $CVCP_3$ problem the constant $c$ can be $4$.

We will solve more general versions of the $CVCP_3$ problem where additionally as a part of the input we are given a set $S\subseteq V$ which contains vertices that must belong to a solution.


%
%

\begin{center}
\begin{tabular}{lp{11cm}}
\multicolumn{2}{l}{\textsc{Constrained Connected Vertex Cover $P_3$} (Constrained CVCP${}_3$)}\\[1ex]
\textbf{Input:} & A graph $G=(V,E)$, a subset $S \subseteq V$ and an integer $k$. \\
\textbf{Question:} & Does there exist a VCP${}_3$-set $F$ of size at most $k$ in $G$ such that $S \subseteq F$ and $G[F]$ is connected?
\end{tabular}
\end{center}

\noindent {\bf Remark.} We assume that the set $S\subseteq V$ in Constrained $CVCP_3$ is nonempty. To solve the problem where $S=\emptyset$, we simply iterate over all possible choices of $v_1\in V$ and put $S=\{v_1\}$.

\begin{defn}
Given a graph $G$, assume $F$ is a $VCP_3$ set of $G$. A cut of $F$ is a pair $(F^1,F^2)$ with $F^1\cap F^2=\emptyset$, $F^1\cup F^2=F$. We refer to $F^1$ and $F^2$ as to the sides of the cut. A cut $(F^1,F^2)$ of $F$ is called a consistent cut of $F$ if $u\in F^1$ and $v\in F^2$ imply $uv\notin E(G[F])$.
\end{defn}

Since no edge of $G[F]$ can go across the cut, each of its connected components has to be either entirely in $F^1$, or entirely in $F^2$. But we can freely choose the side of the cut for every connected component of $G[F]$ independently. Thus, if the induced subgraph $G[F]$ has $c$ connected components, then the $VCP_3$ set is consistent with $2^c$ cuts. Unfortunately if $G[F]$ is connected, $F$ is consistent with two cuts, $(F,\emptyset)$ and $(\emptyset,F)$.
What we do want is that if $G[F]$ is connected, $F$ is consistent with only one cut. Thus, we will select an arbitrary vertex $v_1$ and fix it permanently into the $F^1$ side of the cut.

\begin{thm}
There exists a randomized algorithm that given an $n$-vertex graph $G=(V,E)$ with a tree decomposition of width $p$ solves Constrained $CVCP_3$ in $O(4^p\cdot pn^{6})$ time. The algorithm cannot give false positives and may give false negatives with probability at most 1/2.
\end{thm}

\begin{proof}
We use the Cut\&Count technique \cite{cygan2} which can be applied to problems with certain connectivity requirements. Let $\mathcal{F}\subseteq 2^V$ be a set of
solutions of Constrained $CVCP_3$; we aim to decide whether it is empty. Cut\&Count is split in two parts: (1) The Cut part: Relax the connectivity requirement by considering the set $\mathcal{R}\supseteq \mathcal{F}$ of all $VCP_3$ sets containing $S$. Furthermore, consider the set $\mathcal{C}$ of pairs $(X,C)$ where $X\in \mathcal{R}$ and $C$ is a consistent cut of $F$; (2) The Count part: Compute $|\mathcal{C}|$ modulo 2 using a sub-procedure. Non-connected candidate solutions $X\in \mathcal{R}\setminus \mathcal{F}$ cancel since they are consistent with an even number of cuts. Connected candidates in $\mathcal{F}$ remain.

{\bf The Cut part.}

Suppose we are given a weight function $w\colon V\rightarrow\{1,2,\cdots,N\}$, taking $N=2n$. For a subset $X\subseteq V$, $w(X)=\sum_{v\in X}w(v)$. Recall that we may assume that $S\neq \emptyset$ and we choose one fixed vertex $v_1\in S$. For any integer $W$, we define:
\begin{eqnarray*}
\mathcal{R}_W&=&\{F\subseteq V(G): \text{$F$ is a  $VCP_3$ set of $G$, $|F|=k$, $S\subseteq F$ and $w(F)=W$}\},\\
\mathcal{F}_W&=&\{F\subseteq V(G): \text{$F\in \mathcal{R}_W$ and $G[F]$ is connected}\},\\
\mathcal{C}_W&=&\{(F,(F^1,F^2)): F\in \mathcal{R}_W, (F^1,F^2)\text{ is a consistent cut of $F$, and $v_1\in F^1$}\}.
\end{eqnarray*}

The set $\mathcal{R}_W$ is the set of candidate solutions, where we relax the connectivity requirement. The set $\bigcup_W \mathcal{F}_W$ is our set of solutions. If for any $W$ this set is nonempty, our problem has a positive answer. By the following Lemma we show that instead of calculating the parity of $\mathcal{F}_W$ directly, we can count pairs consisting of a set from $\mathcal{R}_W$ and a cut consistent with it, i.e., $\mathcal{C}_W$.

\begin{lemma}\label{lemma2}
Let $w,\mathcal{F}_W,\mathcal{C}_W$ be as defined above. The parity of $|\mathcal{C}_W|$ is the same as the parity of $|\mathcal{F}_W|$, in other words, $|\mathcal{C}_W|\equiv |\mathcal{F}_W|(\bmod2)$.
\end{lemma}

{\bf Proof of Lemma \ref{lemma2}.}
 For any $VCP_3$ set $F$ from $\mathcal{R}_W$, each connected component of $G[F]$ has to be contained either in $F^1$ or in $F^2$. However, the connected component of $G[F]$ containing $v_1$ has to be contained in the $F^1$ side of the cut. Thus, if $G[F]$ has $c$ components, then $F$ is consistent with $2^{c-1}$ cuts, which is an odd number for $F\in \mathcal{F}_W$ and an even number for $F\in \mathcal{R}_W\setminus \mathcal{F}_W$. \qed\\

\noindent {\bf The Count part.}

Next we describe a procedure $CountC(w,W,\mathcal{T})$ that given a refined variant of nice tree decomposition $\mathcal{T}$, weight function $w$ and an integer $W$, computes $|\mathcal{C}_W|$ modulo 2.

Recall that for each node $t$ of $T$ we denote by $V_t$ the union of all the bags present in the subtree of $T$ rooted at $t$ and

 $$G_t=(V_t,E_t=\{e: e \text{ is introduced in the subtree rooted at t}\}).$$

For every node $t\in T$ of the tree decomposition, integers $0\leq i\leq |V|=n$, $0\leq w\leq N|V|=2n^2$ and colorings (or, mappings) $f\colon X_t\rightarrow\{1_1,1_2,0_0,0_1\}$, define

\begin{eqnarray*}
\mathcal{R}_t(i,w)&=&\{F\subseteq V_t: (\text{F is a $VCP_3$ set of $G_t$})\wedge(|F|=i)\wedge (w(F)=w)\\
&&\wedge (S\cap V_t\subseteq F)\},\\
\mathcal{C}_t(i,w)&=&\{(F,(F^1,F^2)): (F\in \mathcal{R}_t(i,w))\wedge ((F^1,F^2)\text{ is a consistent cut of $F$})\\
&&\wedge(v_1\in V_t\Rightarrow v_1\in F^1)\},\\
A_t(i,w,f)&=&|\{(F,(F^1,F^2))\in \mathcal{C}_t(i,w): (f(v)=1_1\Rightarrow v\in F^1) \wedge (f(v)=1_2\Rightarrow v\in F^2)\\
&&\wedge (f(v)=0_0\Rightarrow\text {$v$ is an isolated vertex in $G_t-F$})\\
&&\wedge (f(v)=0_1\Rightarrow\text {$v$ is a 1-degree vertex in $G_t-F$})\}|.
\end{eqnarray*}

$A_t(i,w,f)$ is the number of pairs from $\mathcal{C}$ of candidate solutions and consistent cuts on $G_t$, with fixed size, weight and interface on vertices from $V_t$. And the number $A_t(i,w,f)$ counts those elements of $\mathcal{C}_t(i,w)$ which additionally behave on vertices of $V_t$ in a fashion prescribed by the coloring $f$. In particular note that

$$\sum_{f}A_t(i,w,f)=|\mathcal{C}_t(i,w)|,$$
where the sum is taken over all colorings from $X_t$ to $\{1_1,1_2,0_0,0_1\}$. For a vertex $x\in X_t$, $f(v)=1_j$ means $v\in F$ and $v\in F^j$ for $j=1,2$, $f(v)=0_0$ means $v$ is an isolated vertex in $G_t-F$, $f(v)=0_1$ means $v$ is a 1-degree vertex in $G_t-F$.

Recall that $r$ is the root of the tree decomposition, we have that $|\mathcal{C}_W|=|\mathcal{C}_r(k,W)|$. Since we are interested in values $\mathcal{C}_W$ modulo 2, it suffices to compute values $A_r(k,W,\emptyset)$ for all $W$.

We present a dynamic programming algorithm that computes $A_t(i,w,f)$ for all nodes $t\in T$ in a bottom-up fashion for all reasonable values of $i$, $w$ and $f$.

\medskip
\noindent {\bf Leaf node.}

For a leaf node $t$, $X_t=\emptyset$. $A_t(0,0,\emptyset)=1$ and all other values of $A_t(i,w,f)$ are zeros.

\medskip
\noindent{\bf Introduce node.}

Let $t$ be an introduce node with a child $t'$ so that $X_t=X_{t'}\cup\{v\}$ for some $v\notin X_{t'}$. Note that $V_t=V_{t'}\cup\{v\}$ and $E_t=E_{t'}$, i.e., $v$ is an isolated vertex in $G_t$.

$$A_t(i,w,f)=\left\{
               \begin{array}{ll}
                 A_{t'}(i-1,w-w(v),f|_{X_{t'}}), & \hbox{when $f(v)=1_1$;} \\
                 A_{t'}(i-1,w-w(v),f|_{X_{t'}}), & \hbox{when $f(v)=1_2$ and $v\neq v_1$;} \\
                 0, & \hbox{when $f(v)=1_2$ and $v=v_1$;} \\
                 A_{t'}(i,w,f|_{X_{t'}}), & \hbox{when $f(v)=0_0$;} \\
                 0, & \hbox{when $f(v)=0_1$,}
               \end{array}
             \right.
$$
where $f|_{X_{t'}}$ is the restriction of $f$ to $X_{t'}$.

\medskip
\noindent{\bf Introduce edge node.}

Let $t$ be an introduce edge node labeled with an edge $uv\in E(G)$ and $t'$ be the child of $t$. Note that $V_t=V_{t'}$, $E_{t}=E_{t'}\cup\{uv\}$ and $uv\notin E_{t'}$.

$$A_t(i,w,f)=\left\{
               \begin{array}{ll}
                 0, & \hbox{when $(f(u),f(v))\in\{(1_1,1_2),(1_2,1_1)\}$;} \\
                 0, & \hbox{when $(f(u),f(v))=(1_2,1_2)$ and $v_1\in\{u,v\}$;} \\
                 0, & \hbox{when $(f(u),f(v))\in\{(0_0,0_0),(0_0,0_1),(0_1,0_0)\}$;} \\
                 A_{t'}(i,w,f'), & \hbox{when $(f(u),f(v))=(0_1,0_1)$;} \\
                 A_{t'}(i,w,f), & \hbox{otherwise,}
               \end{array}
             \right.
$$
where $f'$ is a coloring of $X_t=X_{t'}$, and $f'(w)=\left\{
            \begin{array}{ll}
              f(w), & \hbox{when $w\neq u$ and $w\neq v$;} \\
              0_0, & \hbox{when $w=u$ or $w=v$.}
            \end{array}
          \right.
$

\medskip
\noindent{\bf Forget node.}

Let $t$ be a forget node with a child $t'$ so that $X_t=X_{t'}\setminus\{v\}$ for some $v\in X_{t'}$. Note that $V_t=V_{t'}$ and $G_t=G_{t'}$. For a coloring $f$ of $X_t$ and a color $\alpha\in\{1_1,1_2,0_0,0_1\}$, we define a coloring $f_{v\rightarrow \alpha}$ of $X_{t'}$ as follows:

$$f_{v\rightarrow\alpha}(x)=\left\{
         \begin{array}{ll}
           f(x), & \hbox{when $x\neq v$,} \\
           \alpha, & \hbox{when $x=v$.}
         \end{array}
       \right.
$$

We have
$$A_t(i,w,f)=\sum_{\alpha\in\{1_1,1_2,0_0,0_1\}}A_{t'}(i,w,f_{v\rightarrow\alpha}).$$

\medskip
\noindent{\bf Join node.}

Suppose that $t$ is a join node with children $t_1$, $t_2$ so that $X_t=X_{t_1}=X_{t_2}$. Note that $V_t=V_{t_1}\cup V_{t_2}$, $V_{t_1}\cap V_{t_2}=X_t$ and $E_t=E_{t_1}\cup E_{t_2}$.

We say that a pair of colorings $f_1$ of $X_{t_1}$, $f_2$ of $X_{t_2}$ is consistent with a coloring $f$ of $X_t$ if for each vertex $v\in X_t$, the following conditions hold.
\begin{description}
  \item[(1)] $f(v)=1_1$ if and only if $f_1(v)=f_2(v)=1_1$,
  \item[(2)] $f(v)=1_2$ if and only if $f_1(v)=f_2(v)=1_2$,
  \item[(3)] $f(v)=0_0$ if and only if $f_1(v)=f_2(v)=0_0$,
  \item[(4)] $f(v)=0_{1}$ if and only if $(f_1(v),f_2(v))\in\{(0_{1},0_0),(0_0,0_{1})\}$.
\end{description}

We can now write a recursion formula for join nodes.

$$A_t(i,w,f)=\sum_{i_1+i_2=i+|f^{-1}(\{1_1,1_2\})|}\ \ \sum_{w_1+w_2=w+w(f^{-1}(\{1_1,1_2\}))}\ \sum_{f_1,f_2}A_{t_1}(i_1,w_1,f_1)A_{t_2}(i_2,w_2,f_2),$$
where the third sum is taken over all pairs of colorings $f_1$, $f_2$ consistent with $f$. To achieve the coloring $f$, the colorings of children must be a pair of colorings consistent with $f$.

On one hand, if $(F,(F^1,F^2))\in \mathcal{C}_t(i,w)$ is a pair of candidate solutions for coloring $f$ and consistent cuts on $G_t$, with fixed size $i$, weight $w$ and interface on vertices from $V_t$, then $F_1=F\cap V_{t_1}$ and $F_2=F\cap V_{t_2}$ are $VCP_3$ sets for $t_1$ and $f_1$, and $t_2$ and $f_2$, for some pairs of colorings $f_1,f_2$ that are consistent with $f$. Accordingly, $(F_k,(F^1_k,F^2_k))$ ($k=1,2$) is a pair of candidate solutions and consistent cuts on $G_{t_k}$, with fixed size $i_k$, weight $w_k$ and interface on vertices from $V_{t_k}$, where $F_k^j=\{v|v\in F_k \wedge f_k(v)=1_j\}$ for $j=1,2$.
Since vertices coloured $1_j$ for $j=1,2$ in $X_t$ are accounted for both numbers $A_{t}(i,w,f)$ of the children, we add their contribution to the accumulators.

On the other hand, given a coloring $f$ of $X_t$, consider any pair of colorings $f_1,f_2$ that is consistent with $f$. Suppose $(F_k,(F^1_k,F^2_k))$ ($k=1,2$) is a pair of candidate solutions for coloring $f_k$ and consistent cuts on $G_{t_k}$, with fixed size $i_k$, weight $w_k$ and interface on vertices from $V_{t_k}$. Then $(F_1\cup F_2,(F^1_1\cup F^1_2,F^2_1\cup F^2_2))$ is pair of candidate solutions for coloring $f$ and consistent cuts on $G_t$, with fixed size $i_1+i_2-|f^{-1}(\{1_1,1_2\})|$, weight $w_1+w_2-w(f^{-1}(\{1_1,1_2\}))$ and interface on vertices from $V_t$.

So, the above recursive formula holds.

Similarly, we can show that the above recurrences lead to a dynamic programming algorithm that computes the parity of $|\mathcal{F}_W|$ for all values of $W$ in $O(4^p\cdot pn^{6})$ time, since $|\mathcal{C}_W|=A_r(k,W,\emptyset)$ and $|\mathcal{F}_W|\equiv|\mathcal{C}_W|$.

It remains to show why solving the parity version of the problem is helpful when facing the decision version. We use the following isolation lemma.

\begin{defn}
A function $w\colon U\rightarrow\mathbb{Z}$ isolates a set family $\mathcal{F}\subseteq 2^U$ if there is a unique $S'\in \mathcal{F}$ with $w(S')=\min_{S\in \mathcal{F}}w(S)$.
\end{defn}

\begin{lemma}\cite{cygan,cygan2}
{\bf (Isolation Lemma).} Let $\mathcal{F}\subseteq 2^U$ be a set family over a universe $U$ with $|\mathcal{F}|>0$. For each $u\in U$, choose a weight $w(u)\in\{1,2,\cdots,N\}$ uniformly and independently at random. Then

$$Pr(w \text{ isolates } \mathcal{F})\geq 1-\frac{|U|}{N}.$$
\end{lemma}

The Isolation Lemma gives a reduction from the decision version of the problem to the parity version as follows. Our universe is $U=V(G)$, and we put $$\mathcal{F}=\{F\subseteq V(G):F \text{ is a $VCP_3$ set of $G$, $|F|\leq k$, $S\subseteq F$ and $G[F]$ is connected.}\}.$$ We sample an integral weight $w\colon V\rightarrow\{1,2,\cdots,N=2|U|\}$. For each $W$ between 1 and $2nk$ we compute the parity of the number of pairs in $\mathcal{C}_W$. Since $\mathcal{C}_W\equiv\mathcal{F}_W$, if for some $W$ the parity is odd, we answer YES, otherwise we answer NO.

Even though there might be exponentially many solutions, that is elements of $\mathcal{F}$, with probability at least 1/2 for some weight $W$ between 1 and $2kn$ there is exactly one solution from $\mathcal{F}$ of weight $W$ with respect to $w$. Thus we can obtain a Monte Carlo algorithm where the probability of error is at most 1/2.

We have completed the proof of the theorem.\end{proof}

\acknowledgements

The authors are grateful to the anonymous referees for their comments and constructive suggestions that
greatly improved the original manuscript.

\nocite{*}
\bibliographystyle{abbrvnat}
\bibliography{finalversion}

\begin{thebibliography}{23}
\providecommand{\natexlab}[1]{#1}
\providecommand{\url}[1]{\texttt{#1}}
\expandafter\ifx\csname urlstyle\endcsname\relax
  \providecommand{\doi}[1]{doi: #1}\else
  \providecommand{\doi}{doi: \begingroup \urlstyle{rm}\Url}\fi

\bibitem[Alber and Niedermeier(2002)]{alber2}
J.~Alber and R.~Niedermeier.
\newblock Improved tree decomposition based algorithms for domination-like
  problems.
\newblock In \emph{LATIN}, pages 613--628, 2002.

\bibitem[Bj\"orklund et~al.(2007)Bj\"orklund, Husfeldt, Kaski, and
  Koivisto]{bjo}
A.~Bj\"orklund, T.~Husfeldt, P.~Kaski, and M.~Koivisto.
\newblock Fourier meets m\"{o}bius: Fast subset convolution.
\newblock In \emph{STOC}, pages 67--74, 2007.

\bibitem[Bodlaender(1996)]{bod1}
H.~Bodlaender.
\newblock A linear time algorithm for finding tree-decompostions of small
  treewidth.
\newblock \emph{SIAM J. Comput.}, 25:\penalty0 1305--1317, 1996.

\bibitem[Bodlaender and Koster(2010)]{bod3}
H.~Bodlaender and A.~Koster.
\newblock Treewidth computations i. upper bounds.
\newblock \emph{Inform. Comput.}, 208:\penalty0 259--275, 2010.

\bibitem[Bondy and Murty(2008)]{bon}
J.~Bondy and U.~Murty.
\newblock \emph{Graph Theoryd}.
\newblock Springer, 2008.
\newblock New York.

\bibitem[Chang et~al.(2014)Chang, Chen, Hung, Liu, Rossmanith, and
  Sikdar]{chang}
M.~Chang, L.~Chen, L.~Hung, Y.~Liu, P.~Rossmanith, and S.~Sikdar.
\newblock An $o^*(1.4658^n)$-time exact algorithm for the maximum
  bounded-degree-1 set problem.
\newblock In \emph{CMCT}, volume 244, pages 9--18, 2014.

\bibitem[Cygan et~al.(2011)Cygan, Nederlof, Pilipczuk, Pilipczuk, Rooij, and
  Wojtaszczyk]{cygan2}
M.~Cygan, J.~Nederlof, M.~Pilipczuk, M.~Pilipczuk, J.~Rooij, and
  J.~Wojtaszczyk.
\newblock Solving connectivity problems parameterized by treewidth in single
  exponential time.
\newblock In \emph{FOCS}, pages 150--159, 2011.
\newblock Full version is available as arXiv:1103.0534.

\bibitem[Cygan et~al.(2015)Cygan, Fomin, Kowalik, Lokshtanov, Mark, Pilipczuk,
  Pilipczuk, and Saurabh]{cygan}
M.~Cygan, F.~Fomin, L.~Kowalik, D.~Lokshtanov, D.~Mark, M.~Pilipczuk,
  M.~Pilipczuk, and S.~Saurabh.
\newblock \emph{Parameterized Algorithms}.
\newblock Springer, 2015.
\newblock New York.

\bibitem[Fellows et~al.(2011)Fellows, Guo, Moser, and Niedermeier]{fel}
M.~Fellows, J.~Guo, H.~Moser, and R.~Niedermeier.
\newblock A complexity dichotomy for finding disjoint solutions of vertex
  deletion problems.
\newblock \emph{ACM Trans. Comput. Theory}, 2:\penalty0 1--23, article 5, 2011.

\bibitem[Kardo\v{s} et~al.(2011)Kardo\v{s}, Katreni\v{c}, and Schiermeyer]{kar}
F.~Kardo\v{s}, J.~Katreni\v{c}, and I.~Schiermeyer.
\newblock On computing the minimum 3-path vertex cover and dissociation number
  of graphs.
\newblock \emph{Theoret. Comput. Sci.}, 412:\penalty0 7009--7017, 2011.

\bibitem[Katreni\v{c}(2016)]{kat}
J.~Katreni\v{c}.
\newblock A faster fpt algorithm for 3-path vertex cover.
\newblock \emph{Inform. Process. Lett.}, 116:\penalty0 273--278, 2016.

\bibitem[Kleinberg and Tardos(2005)]{klein}
J.~Kleinberg and E.~Tardos.
\newblock \emph{Algorithm Design}.
\newblock Addison-Wesley, 2005.

\bibitem[Kloks(1994)]{klo}
T.~Kloks.
\newblock Treewidth, computations and approximations.
\newblock In \emph{Lecture Notes in Computer Science}, volume 842. Springer,
  1994.

\bibitem[Lewis and Yannakakis(1980)]{lewis}
J.~Lewis and M.~Yannakakis.
\newblock The node-deletion problem for hereditary properties is np-complete.
\newblock \emph{J. Comput. Syst. Sci.}, 20:\penalty0 219--230, 1980.

\bibitem[Moser et~al.(2012)Moser, Niedermeier, and Sorge]{moser}
H.~Moser, R.~Niedermeier, and M.~Sorge.
\newblock Exact combinatorial algorithms and experiments for finding maximum
  $k$-plexes.
\newblock \emph{J. Combin. Optim.}, 24:\penalty0 347--373, 2012.

\bibitem[Niedermeier(2006)]{nie}
R.~Niedermeier.
\newblock \emph{Invitation to Fixed-Parameter Algorithms}.
\newblock Oxford Univ. Press, 2006.

\bibitem[Pilipczuk(1994)]{pil}
M.~Pilipczuk.
\newblock Problems parameterized by treewidth tractable in single exponential
  time: a logical approach.
\newblock In \emph{MFCS}, pages 520--531, 1994.

\bibitem[Telle and Proskurowski(1993)]{telle}
J.~Telle and A.~Proskurowski.
\newblock Practical algorithms on partial $k$-trees with an application to
  domination-like problemsh.
\newblock In \emph{WADS}, pages 610--621, 1993.

\bibitem[Tu and Zhou(2011{\natexlab{a}})]{tu2}
J.~Tu and W.~Zhou.
\newblock A primal-dual approximation algorithm for the vertex cover $p_3$
  problem.
\newblock \emph{Theoret. Comput. Sci.}, 412:\penalty0 7044--7048,
  2011{\natexlab{a}}.

\bibitem[Tu and Zhou(2011{\natexlab{b}})]{tu3}
J.~Tu and W.~Zhou.
\newblock A factor 2 approximation algorithm for the vertex cover $p_3$
  problem.
\newblock \emph{Inform. Process. Lett.}, 111:\penalty0 683--686,
  2011{\natexlab{b}}.

\bibitem[Tu et~al.(2017)Tu, Wu, Yuan, and Cui]{tu4}
J.~Tu, L.~Wu, J.~Yuan, and L.~Cui.
\newblock On the vertex cover $p_3$ problem parameterized by treewidth.
\newblock \emph{J. Combin. Optim.}, 34:\penalty0 414--425, 2017.

\bibitem[van Rooij et~al.(2009)van Rooij, Bodlaender, and Rossmanithi]{van}
J.~van Rooij, H.~Bodlaender, and P.~Rossmanithi.
\newblock Dynamic programming on tree decompositions using generalised fast
  subset convolution.
\newblock In \emph{ESA}, pages 566--577, 2009.

\bibitem[Wu(2015)]{wu}
B.~Wu.
\newblock A measure and conquer approach for the parameterized bounded
  degree-one vertex deletion.
\newblock In \emph{COCOON}, pages 469--480, 2015.

\end{thebibliography}

\end{document}